\documentclass[a4paper,12pt]{amsart}

\usepackage{amssymb,amsbsy,amsmath,amsfonts,amssymb,amscd}
\usepackage{latexsym}
\usepackage{amsxtra}
\usepackage{amscd}
\usepackage{graphics}
\usepackage{epic}
\usepackage{color}
\input xy
\xyoption{all}

\usepackage{ mathrsfs, amsfonts}

\newcommand\sC{{\mathcal C}}

\newcommand\sD{{\mathcal D}}
\newcommand\sE{{\mathcal E}}

\newcommand\sF{{\mathcal F}}

\newcommand\sL{{\mathcal L}}

\newcommand\sB{{\mathcal B}}

\newcommand\la{\lambda}

\newcommand\Ga{\Gamma}

\newcommand\ga{\gamma}

\newcommand\be{\beta}

\newcommand\de{\delta}

\newcommand{\CC}{\ensuremath{\mathbb{C}}}
\newcommand{\RR}{\ensuremath{\mathbb{R}}}
\newcommand{\ZZ}{\ensuremath{\mathbb{Z}}}
\newcommand{\QQ}{\ensuremath{\mathbb{Q}}}

\newcommand{\NN}{\ensuremath{\mathbb{N}}}
\newcommand{\hol}{\ensuremath{\mathcal{O}}}

\newcommand{\HH}{\ensuremath{\mathbb{H}}}
\newcommand{\PP}{\ensuremath{\mathbb{P}}}

\newcommand{\ra}{\ensuremath{\rightarrow}}

\def\eea{\end{eqnarray*}}
\def\bea{\begin{eqnarray*}}

\newcommand\dual{\mathrel{\raise3pt\hbox{$\underline{\mathrm{\thinspace d
\thinspace}}$}}}
\newcommand\qe{\ifhmode\unskip\nobreak\fi\quad $\Box$}       

\def\BOX{\hfill\lower.5\baselineskip\hbox{$\Box$}}

\newtheorem{theo}[equation]{Theorem}
\newtheorem{remarkk}[equation]{Remark}

\newtheorem{defin}[equation]{Definition}

\newtheorem{lemma}[equation]{Lemma}
\newtheorem{example}[equation]{Example}

\newcommand{\al}{{\alpha}}

\begin{document}

\title[Monodromy and normal forms ]{ Monodromy and normal forms }
\author{ Fabrizio Catanese}
\address {Lehrstuhl Mathematik VIII -Lehrstuhl Mathematik IV\\
Mathematisches Institut der Universit\"at Bayreuth\\
NW II,  Universit\"atsstr. 30\\
95447 Bayreuth}
\email{Fabrizio.Catanese@uni-bayreuth.de}

\thanks{
The present work took place in the realm of the ERC advanced grant TADMICAMT}

\date{\today}

\maketitle

\begin{abstract}
We discuss the history of the monodromy theorem, starting from Weierstra\ss, and the concept of monodromy group. From this viewpoint we compare
then  the Weierstra\ss , the Legendre and other  normal forms for elliptic curves,
explaining their geometric meaning and distinguishing them by their stabilizer in $\PP SL(2,\ZZ)$ and their monodromy.
Then we focus on the birth of the concept of the Jacobian variety, and the geometrization of the theory of Abelian functions and integrals.
We end illustrating the methods of complex analysis  in the simplest issue, the difference equation $f(z) = g(z+1) - g(z)$ on $\CC$.

\end{abstract}

\tableofcontents
\section*{Introduction}
In Jules Verne's novel of 1874, `Le Tour du monde en quatre-vingts jours' , Phileas Fogg is led to his remarkable adventure 
by a bet made in his Club: is it possible to make a tour of the world in 80 days?

Idle questions and bets can be very stimulating, but very difficult to answer when they deal with the history of mathematics,
and one asks how certain ideas, which have been a common knowledge for long time, did indeed evolve and mature through a long period of time,
and through the contributions of many people.

In short, there are  three idle questions which occupy my attention since some time:

\begin{enumerate}
\item
When was the statement of the monodromy theorem first fully formulated (resp. : proven)?
\item
When did  the normal form for elliptic curves
$$ y^2 = x (x-1)(x - \la),$$ which is by nowadays'  tradition called by many (erroneously?)  `the Legendre normal form' first appear?
\item
The old `Jacobi inversion theorem' is today geometrically formulated through the geometry of the  `Jacobian variety $J(C)$' of an algebraic curve $C$ of genus $g$:
when did this formulation clearly show up (and so clearly that, ever since, everybody was talking only in terms of the Jacobian variety)?

\end{enumerate}

The above questions not only deal with themes of research which  were central to Weierstra\ss' work on complex function theory,
but indeed they single out philosophically the importance in mathematics of clean formulations and rigorous arguments.

Ath his point it seems appropriate to cite Caratheodory, who wrote so in the preface of his two volumes on `Funktionentheorie' (\cite{caratheodory}):

` The genius of B. Riemann (1826-1865) intervened not only to bring the Cauchy theory to a certain completion, but also to create the foundations for
 the geometric theory of functions. At almost the same time, K. Weierstra\ss (1815-1897) took up again the above-mentioned idea of Lagrange's
 \footnote{`whose bold idea was to develop the entire theory on the basis of power series'; the definition of an {\bf analytic} function as one which is locally the sum of a power series is due to Lagrange in \cite{lagrange},  please observe that Weierstra\ss \ started to develop the method of power series quite early,
 around 1841 }, on the basis of which he was able to arithmetize 
 Function Theory and to develop a system that in point of rigor and beauty cannot be excelled. The Weierstra\ss \  tradition was carried on in an especially pure form by A. Pringsheim (1850-1941),
 whose book (1925-1932) is extremely instructive.'
 
 Then Caratheodory comments first on the antithesis:
 
 `During the last third of the 19th Century the followers of Riemann and those of Weierstra\ss \  formed two sharply separated schools of thought.'\footnote{and things 
 were made more complicated by some sort of  direct rivalry between the G\"ottingen and Berlin schools of mathematics.}
 
 and then on the sinthesis:
 ` However,
 in the 1870's Georg Cantor (1845-1918) created the Theory of Sets. .. With the aid of Set Theory it was possible for the concepts and results of Cauchy's and Riemann's theories to be put on just as firm basis as that on which Weierstra\ss \ ' theory rests, and this led to the discovery of great new results in the Theory of Functions as well as of many simplifications in the exposition.'
 
 Needless to say, the great appeal of Function Theory rests on two aspects: the fact that classical functions of a real variable are truly understood only after
 one extends their definition to the complex domain, seeing there the maximal domain where the function extends without acquiring singularities.
 
 And, even more, the variety of different methods and perspectives: polynomials, power series, analysis, and geometry, all of  these illustrate several facets
 of the theory of holomorphic functions, complex differentiability, analiticity (local representation through a power series), conformality.
 
 As a concrete example of the several souls which are indispensable  in order to treat  problems in Function Theory, I shall then illustrate a simple case of a crucial
 technical result in the  theory of periodic functions, looking at the different methods which can be used:
 
 given a holomorphic function $f(z)$ on the entire complex plane $\CC$, does there exist a function $g(z)$ solution of the following difference equation?
 $$ (*) \ g(z+1) - g(z) = f(z) .$$

\section{The monodromy theorem}

Curiously enough, there are several famous monodromy theorems, the classical one and some  modern ones.

The classical one is easier to understand, it revolves around the concept of function, and distinguishes between a `monodromic' (i.e., `single valued')
function, and a `polydromic' (or `multiple valued') function.

Typical examples of monodromic functions are the `functions' in modern (Cantor's)  sense:
$$ z \mapsto z^2 , \  z \mapsto e^z : = \Sigma_n \frac{z^n}{n!}, \  z \mapsto cos (z) : = \frac{1}{2} (e^{iz} + e^{-iz}),$$
examples of polydromic functions are
$$ z \mapsto \sqrt{z}, z \mapsto \sqrt{(1-z^2)(1-k^2 z^2)},$$
$$ \ z \mapsto log(z) : = \int_1^z \frac{dt}{t},  \ z \mapsto arcsin (z) : = \int_0^z  \frac{dt}{\sqrt{(1-t^2)}}$$
$$z \mapsto \sL^{-1} (z) : = \int_0^z \frac{dt}{\sqrt{(1-t^2)(1-k^2 t^2)}},$$
where the last function is the Legendre elliptic integral.

The meaning of the two words is easily understood if we recall that in Greek $\mu o \nu o \sigma$ means `single', $\pi o \la y$ means `many',
and $\de \rho o \mu \epsilon \iota \nu$ means `to run'.
So, a function is polydromic if running around some closed path we end up with a  value different from  the beginning one.

In the first two examples, we have an {\bf algebraic function} $y$, which means that there exists a polynomial equation satisfied 
by $z$ and $y$, such as, respectively: 
$$  y^2 = z, \   y^2 = (1-z^2)(1-k^2 z^2),$$
and indeed in modern mathematics we view the algebraic  function as defined by  the second coordinate projection on the curve
$$ f : C : =  \{ (z,y) | y^2 = (1-z^2)(1-k^2 z^2)\} \ra \CC , f(z,y) := y.$$
Here, the reason for the multivaluedness of an algebraic function $y$, defined by an irreducible polynomial equation $P(z,y) = 0$,
is explained more or less by algebra; since, if the degree of $P$ with respect to the variable $y$ equals $d$, then to each value
of $z$ which is not a root of the discriminant $\de(z) = disc_y(P)(z)$ of $p$ with respect to the variable $y$,
there correspond exactly $d$ values of $y$.

The latter examples are deeper to deal with, their nature is transcendental, but in some sense are easier to understand:
one sees indeed that going around the origin the value of the logarithm gets changed by a multiple of
$$ 2 \pi i = \int_{|z| = r} \frac{dt}{t}$$
for each tour $ z = r e^{ i \theta}, 0 \leq \theta \leq 2 \pi$, around the origin; hence the function is polydromic.

The monodromy theorem is deeply based on the concept of {\bf analytic continuation} introduced by Weierstra\ss \  in his lectures (\cite{hurwitz-w}, chapter 10,
pages 93-97). Weierstra\ss \  observes first that  a power series defines inside its convergence disk $D$ a function which is analytic, i.e., it can be represented,
for each point $c \in D$, as a local power series $\Sigma_n (a_n (z-c)^n$; and then that, if two power series $f : D \ra \CC$, $g: D' \ra \CC$
yield the same local power series at some point $c \in D \cap D'$, then the same holds for each other point $c'  \in D \cap D'$. Proceeding in this way,
one can, given a path $\ga$ and a series of disks $D_1, \dots D_r$ such that
$$  D_i \cap D_{i+1} \neq \emptyset, \ \ga \subset \sD : = \cup _i^r D_i ,$$
define the analytic continuation of $f: = f_1$ along the path if there are power series $f_i$ on each disk $D_i$,
such that there is a point $c_i \in D_i \cap D_{i+1}$ where the local power series developments coincide. 

Writes then Hurwitz in his notes of Weierstra\ss \ ' lectures  (page 97):

`L\"a\ss t sich f\"ur einen Punkt nur ein einziges Funktionenelement aufstellen, so hei\ss t die Funktion eindeutig, in entgegengesetzten Falle mehrdeutig.'

The monodromy theorem gives a sufficient condition for the analytic continuation to be single valued (monodromic).

\begin{theo}{(\bf  Monodromy Theorem, Der Monodromie Satz.)}
Let two continuous paths $\ga(s), 0 \leq 1$ and $\de(s), 0 \leq 1$ be given, which have the same end points $\ga(0) = \de(0), \ga(1) = \de(1)$
and which are homotopic, i.e., there is a continuous deformation of $\ga$ and $\de$ given by a continuous function
$ F(s,t) ,\ 0 \leq  s,t \leq 1$ such that $F(s,0) = \ga(s), \ F(s,1) = \de(s)$.  Then analytic continuation along $\ga$ yields the same result as analytic continuation along $\de$. 

In particular, if we have a function $f$ which admits analytic continuation over the whole domain $\Omega$, and the domain $\Omega$, is simply connected,
then $f$ extends to a monodromic (single valued) function.
\end{theo}

Eberhard Freitag states clearly  in his book (\cite{freitag}) that the monodromy theorem was first proven in the lectures of Weierstra\ss \ . 
Indeed, in section 17.2 (pages 136-138 ibidem) it is proven that analytic continuation along segments and triangles gives a single valued result,
and something more is said in chapter 18 (pages 143-145), where it is written that `every analytic function can be made single valued after 
restrriction of its domain of definition': this  amounts, in the explanation by Peter Ullrich (page xxii of the foreword), to showing that, if we take the
star of a point $p$ (union of all segments $[p,z]$ contained in the domain of definition), there any function is single valued.
According to Ullrich, the full statement of the Monodromy theorem for simple connected domains is contained in the `Mitschrift' of Killing (notes of
the Weierstra\ss \  lectures in the summer term 1868), even if the proof is not given in a complete way.

However, the monodromy  theorem does not appear either   in Hurwitz   (see \cite{hurwitz}, which is without the addition by Courant),
nor in Osgood's treatise (\cite{osgood}, and not even in Hermann Weyl's `Die Idee der Riemanschen Fl\"ache' (\cite{weyl}.

The first book source I was able to find the full statement and a complete proof is the book   \cite{h-c}, which appeared in 1922 (already in 1925 there was a second edition!).
The book, which is entitled `Funktionentheorie' has become and remained a classic, it  is the third volume of the Springer series `Grundlehren',
and the monodromy theorem appears (page 348, in chapter 5, entitled `Analytic continuation and Riemann surfaces')
in the section `Geometrische Funktionentheorie' added by Richard Courant (editor of the book) to the lectures by Adolf Hurwitz
`Vorlesungen \"uber allgemeine Funktionentheorie und elliptische Funktionen'.

It is clear that there are two main ingredients in the monodromy theorem: first, the concept of analytic continuation,
second the topological idea of a simply connected domain. It is no coincidence that Courant formulates the theorem for
Riemann surfaces: in our opinion the monodromy theorem represents the ideal marriage of the ideas of Weierstra\ss \  with those
of Riemann.

\subsection{Riemann domain and sheaves}
The concept of analytic continuation is nowadays fully understood through the theory of sheaves, and especially the concept of 
the `espace \'etale' of a sheaf \footnote{`schlicht Gebiet' in German, jokingly called  in Italian by one of my teachers: `spazio lasagnato'.}. 

Sheaves were invented by Leray during the second world war as a way to analyse the topological obstructions to determine global solvability
once local solvability is no problem. For example, any differential form $\phi$ on a differentiable manifold $M$ which is closed (i.e., $d \phi=0$) is locally
exact (so called Poincar\'e' lemma), i.e., locally there exist another form $\psi$ such that $d \psi = \phi$. The De Rham theory shows that in order to 
obtain a globally defined form  $\psi$ such that $d \psi = \phi$ it is necessary and sufficient that the cohomology class of $\phi$ is zero:
this means that the necessary condition $$\int_{N} \phi = 0, \forall N, s.t. \ \partial N = \emptyset $$ is also sufficient.

Today even sheaves are considered by some mathematicians as  elementary mathematics, and the new trend is to disregard
the thorough treatment done in the 1950's and 1960's (\cite{fac}, \cite{godement}) preferring the more general concept of
derived categories (introduced by Grothendieck and Verdier, see  \cite{verdier}).

However, the theory of sheaves associates to a topological space $X$ another topological space $\sF$, with a continuous
map $p$ onto $X$, which is  a local homeomorphism. The points of $\sF$ lying over a point $x \in X$ 
are germs of functions around $x$, which means equivalence classes of functions defined in some neighbourhood $U$ of $x$;
where two germs $f_x, g_x$ are equivalent if there is a smaller neighbourhood $V$ of $x$ where the functions do coincide identically.

This applies in particular to holomorphic functions, which enjoy the special property that the locus where they coincide
identically is both open and closed, and analytic continuation of a holomorphic function $ f : U \ra \CC$ defined over a connected open set $U \subset X$
( $X$ is here a complex manifold) means that we consider first  the open set $U_f$ in $\sF$ given by the germs $f_x$, for $x \in U$, and then we take a larger connected
open set $W \subset \sF$.   The largest such open set (the connected component of $\sF$ containing $U_f$ is called the
{\bf Riemann domain} $\tilde{U}_f$ of the function $f$: it is a complex manifold endowed with a locally biholomorphic map $ p : \tilde{U}_f \ra X$,
an explicit biholomorphism $U \cong U_f$, and a holomorphic function $F : \tilde{U}_f \ra \CC$, which extends the holomorphic function 
determined by $f$ on $U_f$.

\subsection{Monodromy or polydromy?}
As mentioned earlier, there are many monodromy theorems: but most of them deal with a concept, monodromy,
which should be instead called polydromy.  In fact, in the modern use, when we talk about the monodromy of a covering
space, say of an algebraic function, we talk about describing its polydromy, that is, its lack of monodromy.

All this started quite early: see as example for this use page 666 of the encyclopedia by  Wirtinger and others (\cite{wirtinger},
where the `monodromy' of the branch points is being considered), chapter 81 of Bianchi's book \cite{bianchi},
devoted to the monodromy group and the results of Hurwitz (\cite{hurmoves}, where the second part is entitled: Monodromy groups;\footnote{ 
 Fricke and Klein instead distinguish themselves by using the word `polymorphic' for functions which are
 polydromic (\cite{f-k}, vol. II, page 43).}
Frans Oort gave me the following explanation: the monodromy theorem had just become so famous that everybody
felt it so great to talk about it, so that even in the cases where there was polydromy, they fell into the habit
of only talking about monodromy. I have to admit that even nowadays there are slogans which take over so much
in the imagination of mathematicians, that many  would like to understand everything through a few slogans,
forgetting about the complexity of the mathematical world. \footnote{Carl Ludwig Siegel wrote once to Weil: 
`It is completely clear to me which conditions caused the gradual decadence of mathematics, from its high level some 100 years ago, down to the present hopeless nadir... Through the influence of textbooks like those of Hasse, Schreier and van der Waerden, the new generation was seriously harmed, and the work of Bourbaki finally dealt the fatal blow.'  What would Siegel write today? }

For example, the  theory of covering spaces was invented to 
clarify the concept of an algebraic
function and its polydromy.

In the modern terminology one can describe an 
algebraic function $f$ on an algebraic curve $Y$ (equivalently, $Y$ is a compact Riemann surface) 
as a rational
function
$f$ on a projective curve $X$ which admits a 
holomorphic map $p : X \ra Y$ (and we may assume  that $f$ 
generates the corresponding  extension for the respective fields of meromorphic functions
$ \CC(Y) \subset \CC(X)$.

The easiest example would be the one where $Y 
=\PP^1 = \PP^1_{\CC}$ and $f= \sqrt P(x)$,
    $P$ being a square free polynomial.

$f$ is in general polydromic, i.e., many valued as a function on $Y$, and
going around a closed loop we do not return to the same value.
It is a theorem of Weierstra\ss \  (later generalized by Hurwitz in \cite{hur-rat}) that $f$ is a rational
function on $Y$ if $f$ is monodromic, i.e., there is no polydromy.

We shall also adhere here to the prevailing attitude, to call  monodromy what 
should really be called polydromy, and we explain it now in the particular example of algebraic functions
(observe that even in the book by Bianchi \cite{bianchi}, where the author pays special attention to distinguish monodromic and polydromic functions,
as in page 254, 
chapter 81 is entitled: `Gruppo di monodromia').

Given $p : X \ra Y$ as above there is a finite 
set $\sB \subset Y$, called the branch locus,
such that, setting  $Y^* : = Y \setminus \sB, X^* 
: = p^{-1} (Y^* )$, then $p$ induces a covering
space $X^* \ra Y^*$ which is classified by its 
monodromy $\mu$.
  $\mu$ is a homomomorphism
of the fundamental group of 
$Y^*$, $\pi_1 (Y^*, y_0)$  into the group
of permutations of the fibre $  p^{-1} (y_0)$($\pi_1 (Y^*, y_0)$ is the group introduced by Poincar\'e, whose elements are homotopy classes of closed paths
beginning and ending in a fixed point $y_0$) 

If $X$ is irreducible, and $d$ is the degree of 
$p$, then the image of the monodromy is a 
transitive
subgroup of $\mathfrak S_d$, and
conversely Riemann's existence theorem asserts that for any homomorphism
$ \mu \colon  \pi_1 (Y^*, y_0) \ra \mathfrak S_d$ 
with transitive image we obtain an algebraic 
function
on $Y$ with branch set contained in $\sB$.

  Indeed one can factor the monodromy $ \mu 
\colon  \pi_1 (Y^*, y_0) \ra \mathfrak S_d$
through a surjection to a finite group $G$ 
followed by a permutation representation of $G$, 
i.e., an injective
homomorphism $ G \ra \mathfrak S_d$ with transitive image.  A concrete calculus to determine
monodromies explicitly was developed by Hurwitz in his fundamental work \cite{hurmoves}.

While one of the most striking results was obtained in \cite{Schwarz} by Hermann Amandus Schwarz (another  student of Karl Weierstra\ss \ )
who determined all the cases where the monodromy of the Gauss hypergeometric function is finite
(hence algebraicity follows).

\section{Normal forms and monodromy}

The name of Weierstra\ss \  will ever remain present  in complex variables  and function theory  through   terms like `Weierstra\ss \  infinite products'
 (see \cite{werke}, vol. 2, articles 8 and 11), `Weierstra\ss \ ' preparation theorem'
( see \cite{werke}, vol. 2, page 135, the article `Einige auf die Theorie der analytischen Functionen mehrerer ver\"anderlichen sich beziehende S\"atze), which paved the way to the theory of several complex variables  (see \cite{siegel} for a classical exposition, starting from the preparation theorem
and the unique factorization of holomorphic functions of several variables).

But also through a simple notation, for the Weierstra\ss \ ' $\wp $ function. 

Let us recall the by now standard  definition, which is contained nowadays in every textbook of complex analysis \footnote{even if the original definition, see \cite{werke}, vol. 5, was based on the differential equation satisfied by it, and the  relation with the $\sigma$-function.}.

Let $\Omega$ be any discrete subgroup of $\CC$, generated (i.e., $\Omega = \ZZ \omega_1 +  \ZZ \omega_2$) by a real basis $\omega_1, \omega_2$ of $\CC$. Then Weierstra\ss \  defined his $\wp $ function as:
$$ \wp (z) : =  \frac{1}{z^2} + \Sigma'  (\frac{1}{(z - \omega)^2} - \frac{1}{\omega^2}),$$
where $\Sigma' $ denotes summation over the elements $\omega \in \Omega \setminus \{0\}$.

Moreover, he proved that the meromorphic map given by $$ (\wp (z) :  \wp '(z): 1)$$ gives an isomorphism of the 
elliptic curve $ E_{\Omega} : = \CC / \Omega$  with the curve of degree three in the projective plane,   union
of the affine curve in $\CC^2$ defined by the equation 
$$ (*) \ \{ (x,y) |  y^2 = 4 x^3 - g_2 x - g_3 \},$$ 
with the point at infinity with homogeneous coordinates $(1: 0 : 0)$ (the image of the origin in $\CC^2$).
This point is the neutral element $0_E$ for the group law, where three points are collinear if and only if their sum is equal to $0_E$.

The equation
$$ y^2 = 4 x^3 - g_2 x - g_3$$
is called the Weierstra\ss \  normal form of the elliptic curve. Here, an elliptic curve is a group, induced by addition on $\CC$, in particular 
the datum of an elliptic curve yields  a pair $(E,0_E)$.  The curve is smooth if the polynomial $4 x^3 - g_2 x - g_3$
has three distinct roots, equivalently, its discriminant $ \Delta : = g_2^3 - 27 g_3^2 \neq 0$ (if $\Delta = 0$ we get a curve with a node, whose normalization is
isomorphic to $\PP^1$).

Moreover, up to a homotethy in $\CC$, we can assume that $\Omega = \ZZ \oplus \ZZ \tau$,
where $\tau$ is a point in the upper half plane $\{ \tau \in \CC|  Im (\tau) > 0\}$.

The functions  $g_2 (\omega_1,\omega_2), g_3 (\omega_1,\omega_2)$ (one writes $g_2 (\tau), g_3 (\tau)$ when one restricts to pairs of the form 
 $(\omega_1,\omega_2)= (1, \tau)$) are the fundamental examples of Eisenstein series (\cite{eisenstein}, 7. IV, pages 213-335, published in Crelle, before 1847),
$$  g_2 = 60 \cdot   \Sigma' \frac{1}{\omega^4},   g_3 = 140 \cdot \Sigma'  \frac{1}{\omega^6},$$
as one  learns as an undergraduate  student (I read it in the book by Henri Cartan, \cite{cartan}).

It is clear that a homothety, multiplying each $\omega \in \Omega$ by $c$, has the effect of multiplying $g_m$ by $c^{-2m}$, and that a change of basis in $\Omega$
has no effect whatsoever; this implies that $g_2, g_3$ are automorphic forms for  the group of oriented changes of basis in 
$\Omega$: $\PP SL (2 , \ZZ) = SL (2 , \ZZ)/ \pm I$.

This group acts on the upper half plane 
$$  \phi : = \left( \begin{array}{cc} a & b \\ c & d \end{array}  \right)   \in \PP SL (2 , \ZZ) $$
 by 
$$\phi (\tau) =   \frac{ a  \tau + b} { c \tau + d},$$
and then $$g_{m} (\tau ) = g_{m} (\frac{ a  \tau + b} { c \tau + d}) ( c \tau + d)^{-2m}= g_{m} ( \phi (\tau) ) ( c \tau + d)^{-2m}.$$
This transformation formula says that $g_m$ are automorphic forms  of weight $m$, and indeed $g_2, g_3$ are the most important ones,
since the graded ring of automorphic forms for $\PP SL (2 , \ZZ)$ is the polynomial ring $\CC [ g_2, g_3]$.

In the Weierstra\ss \  normal form the origin  of the elliptic curve is fixed (the neutral element for the addition law). It corresponds to the origin in $\CC$, and to the point at infinity 
(i.e., the unique point with $w=0$) of the projective curve
of homogeneous equation
$$ y^2 w = 4 x^3 - g_2 (\tau) x w^2  - g_3 (\tau) w^3.$$
The three intersection points with the $x$-axis correspond to the three half periods $\frac{1}{2},  \frac{\tau}{2} \frac{\tau+1}{2}$:
these are exactly the points $P$ on the elliptic curve which are  2-torsion points ( this means $ 2 P \equiv  0 (mod \ \Omega)$) and different from the origin.
Observe that the 2-torsion points form a subgroup of $E$, isomorphic to $(\ZZ/2)^2$.

The  Weierstra\ss \  normal form is quite elegant, moreover it shows that the elliptic curves are just the
nonsingular plane cubics, the projective plane curves of degree 3 (since we have a field, \CC,  where $2,3$ are invertible).
It has had therefore a profound influence in number theory, where one prefers the `normal form'
$$ y^2 = x (x-1)(x - \la),  \la \neq 0,1,$$
in which the three roots are numbered, and brought, via a unique affine transformation of $\CC$, to be equal to $0,1,\la$.
 This normal form is essentially explained in \cite{werke}, vol. 6, page 136; in chapter 13, which  dedicated to the degree two transformation leading to the Legendre normal form: here fails
only the letter $\la$ for the transform of the fourth root $a_4$ of a degree $4$ polynomial $R(x)$.

The real issue is that the   `Weierstra\ss \  normal form' is not really a normal form!   In the sense that two elliptic curves which are isomorphic
do not have the same invariants $g_2, g_3$: these indeed, as we saw,  change if one represents $E$ by another point $\tau$ in the upper half plane,
say $ \frac{ a  \tau + b} { c \tau + d}$.

  Indeed, two elliptic curves are, as well known, isomorphic if and only if they have the same $j$-invariant.
  The $j$-function can be calculated in terms of the above forms (notice that the right hand side is automorphic of weight
  $ 6 - 6 = 0$, hence it is a well defined function of $\tau$). 
  
    $$ (**) j (\la) = \frac{4}{27} \frac{(\la^2 - \la + 1)^3}{\la^2 (\la-1)^2} = \frac {g_2^3}{g_2^3 - 27 g_3^2}.$$

The normal form used  by Legendre was instead the normal affine form:
 $$ y^2 = ( x ^2-1) ( x ^2-a^2), a \neq 1, -1.$$ 
 
  Weierstra\ss \ explains clearly in his lectures (see \cite{werke}, vol. 5, especially page 319)  how to obtain form his normal form  the Legendre normal form, once one has found the roots
 of the polynomial $P(x) : =  4 x^3 - g_2 x - g_3$. The same is also lucidly done in the book of  Bianchi  \cite{bianchi}, pages 433-436, for any polynomial $P(x)$ of degree $3$ or $4$).
 
 As the Weierstra\ss \  normal form describes a plane curve isomorphic to the elliptic curve $ E = \CC / \Omega$, image through the 
 Weierstra\ss \ ' $\wp $ function and its derivative, something similar 
 happens with the  Legendre normal form: there exists 
  $\mathcal{L}$, 
 a Legendre function for $E$: $\mathcal{L} \colon E \rightarrow \mathbb{P}^1$, a meromorphic function  
which makes $E$  a double cover of $\PP^1$ branched over the  four distinct ordered points:
$\pm 1, \pm a\in \mathbb{P}^1\setminus \{0, \infty\}$. 

As explained by Weierstra\ss \ , and later in the book by Tricomi \cite{tricomiell}, the Legendre normal form is
very important for the applications of elliptic functions: and essentially because the coefficient $a$
is a well defined function of $\tau$, which has finite monodromy on the space  of isomorphicm classes
of elliptic curves, parametrized by the invariant $ j \in \CC$.  Indeed Weierstra\ss \  in \cite{werke}, vol. 5, devotes a lot of efforts to show how one can pass from
one normal form to the other.

Our point of view here  is to relate these three normal forms (and a fourth one)  with the monodromy point of view, in the sense of Fricke and Klein (\cite{f-k}),
and as well explained by Bianchi (who was a student of Klein) in \cite{bianchi}; i.e., to see that certain functions of $\tau$ describe  a quotient
of the upper half plane by a well determined  subgroup of $\PP SL(2, \ZZ)$.

Our interest is not merely historical: the Legendre normal form has played an important role in some discovery of new algebraic surfaces done by Inoue 
in \cite{inoue}, and in our joint work with Bauer and Frapporti ( see for instance \cite{BC11burniat1} and \cite{bcf}).

Before we proceed, let us observe that the first fundamental theorem about projectivities states that, given distinct points $P_1, P_2, P_3,P_4$ in $\PP^1$,
there exists a unique projective transformation sending $P_1 \mapsto \infty, P_2 \mapsto 0, P_3 \mapsto 1$; and the image of the fourth point 
$P_4$ is the cross ratio $$ \frac{P_4- P_2}{P_4 - P_1} \cdot \frac{P_3- P_1}{P_3 - P_2}.$$

In this terms, we obtain that 
 $\la$ is the cross-ratio of the four points $\mathfrak p (0), \mathfrak p (\frac{1}{2}), \mathfrak p (\frac{\tau}{2}), \mathfrak p (\frac{1 + \tau}{2}), $
 where $\mathfrak p $ is the Weierstra\ss \  function.
 
And the   $j$-invariant is just the only invariant for a group of four distinct points in the projective line (the cross-ratio is well defined
for an ordered fourtuple, but permuting the four points has the effect of exchanging $\la$ with six  values, $\la , \frac{1}{\la}, 1 - \la,
\frac{1}{1-\la},  1 -   \frac{1}{\la},  \frac{\la}{\la -1}$).

  The above equations $(**)$ show easily that $j (\la) = 0 \Leftrightarrow \la^3 + 1 = 0 , \la \neq 1$
   $j (\la) = 1 \Leftrightarrow \la = - 1$. Indeed the holomorphic map $\tau \ra j(\tau)$  is the quotient map
   of the upper half plane for the action of the  group $\PP SL(2, \ZZ)$, and over $\CC \setminus \{0,1\}$
   it is a covering space (this leads, see for instance  \cite{bianchi}, pages 359-362, or  \cite{rudin}, pages 324 and foll.,  to a quick proof of Picard's first theorem that any entire function
   which omits two values is necessarily constant).
   
   Let us now consider the Legendre normal form $y^2 = (x^2-1) (x^2- a^2)$: since this is an elliptic curve $E$, there exists a $\tau$ in the upper half plane,
    such that $E = E_{\tau}$,
   and a function   $\mathcal{L}\colon E \rightarrow \PP^1$ which is called  a {\bf Legendre function} for $E$.
   
   This function is a close relative of the Weierstra\ss \  function, and we can  give this function in dependence of the parameter $a$; in order however to
   see the symmetries of $E$ related to the subgroup of 2-torsion points of $E$,  it is convenient to set $ a = b^2$,
   i.e., to take a square root of $a$.

The Inoue normal form is  the normal form
 $$ y^2 = (\xi^2-1) (\xi^2-b^4).$$ 

\noindent
We have the following relations for the Legendre function (see \cite[Lemma 3-2]{inoue}, and \cite[Section 1]{BC11burniat1} for an algebraic treatment):
\begin{itemize}
\item $\mathcal{L}(0)=1$, $\mathcal{L}(\frac{1}{2})=-1$,  $\mathcal{L}(\frac{\tau}{2})=a$,
 $\mathcal{L}(\frac{\tau+1}{2})=-a$;
\item set $b:=\mathcal{L}(\frac{\tau}{4})$: then $b^2=a$;
\item $\frac{\mathrm d\mathcal L}{\mathrm d z}(z)=0$ if and only if 
$z\in\{0,\frac 12, \frac {\tau}2, \frac {\tau+1}2\}$
 since these are the ramification points of $\mathcal{L}$.
\end{itemize}

\noindent
Moreover,
\begin{eqnarray}
&\mathcal{L}(z)=\mathcal{L}(z+1)=\mathcal{L}(z+\tau)=\mathcal{L}(-z)=
-\mathcal{L}\bigg(z+\dfrac 12\bigg),\nonumber \\
&\mathcal{L}\bigg(z+\dfrac{\tau}2\bigg)= \dfrac{a} {\mathcal L(z)}\,.\nonumber
\end{eqnarray}

The importance of this normal form is to describe explicitly, on the given family of elliptic curves, the action of the group
$(\ZZ/2)^3$ acting by sending $$ z \mapsto \pm z + \frac{1}{2} \omega,  \ \omega \in \Omega / 2 \Omega \cong (\ZZ/2)^2.$$

Observe, from the symmetry point of view, that the Weierstra\ss \  normal form clearly exhibits the symmetry $ (x,y) \mapsto (x,-y)$:
this symmetry corresponds to multiplication by $-1$ on the elliptic curve (i.e., sending a point to its inverse).

On the elliptic curve in Legendre normal form $y^2 = (x^2-1) (x^2 - a^2)$ we have the group $(\ZZ/2)^2$ of automorphisms consisting of 
 $$  g_1 (x,y) = (-x, -y), \ g_2 (x,y) = (x, -y) \  g_3 (x,y) = (-x, y).$$ 
 The transformation $g_1$ corresponds to a translation by a point of 2-torsion, and indeed  
the quotient of $E$  by $g_1$ is easily seen to be the elliptic curve of equation
$$ v^2 = u (u-1) (u - a^2),  \ u : = x^2 ,  \ v : = xy.$$  

One can see  more symmetry via algebra,  considering the following 1-parameter family   of intersections of two quadrics in the projective space $\PP^3$
with homogeneous coordinates $(x_0 : x_1 : x_2 : x_3)$.
 
  $$E  (b) :=\{ x_1^2+x_2^2+x_3^2=0,
 \quad x_0^2= (b^2+1)^2 x_1^2+( b^2-1)^2x_2^2\}\\,$$
 where $ b \in \CC \setminus \{ 0, 1, -1 , i, - i \}$.
 
 On it the group $(\ZZ/2)^3$ acts in a quite simple way, multiplying each variable $x_i$ by $\pm 1$ (we get  the group $(\ZZ/2)^3$ and not  $(\ZZ/2)^4$, 
 since these are projective coordinates, hence $ (-x_0: - x_1 : -x_2 : -x_3) = (x_0: x_1 : x_2 : x_3)$).
 
 The relation with the Legendre normal form is obtained as follows.
 We set $$ (s : t) : = ( x_1 + i x_2 : x_3) = ( - x_3 : x_1 - i x_2), \ \  \xi : = \frac{bs}{t}\,,$$
 and in this way  the family of genus one curves $E(b)$ is the Legendre family of
 elliptic curves in Legendre normal affine form:
 $$ y^2 = (\xi^2-1) (\xi^2-a^2), \ \  a : = b^2\,.$$ 

 \noindent  In fact (see \cite{bcf}):
\begin{eqnarray*} 
  x_0^2&=& (b^2+1)^2 x_1^2+( b^2-1)^2x_2^2 =- (a+1)^2 (s^2-t^2)^2+( a -1)^2(s^2 + t^2) ^2 = \\
 &=& 4 [(a^2 + 1) s^2 t^2  - a (t^4 + s^4)]   = 4 t^4 \left[  (a^2 + 1)\left( \frac{\xi} {b}\right)^2 - a \left( 1 + \left( \frac{\xi} {b}\right)^4\right)\right ] = \\
 & =& - 4 t^4 \frac{1}{b^2} [ - (a^2 + 1) \xi^2  + (a^2 + \xi^4) =  \frac{- 4 t^4}{b^2} [   (\xi^2-1) (\xi^2-a^2)] 
 \end{eqnarray*}
 and it suffices to set $$ y : =  \frac{i b x_0}{2 t^2}\,.$$ 
 
 The group $(\ZZ/ 2)^3$ acts fibrewise on the family $E(b)$ via the commuting involutions:
 $$ x_0 \longleftrightarrow - x_0,  \ \ x_3 \longleftrightarrow  - x_3,\ \  x_1 \longleftrightarrow  - x_1, $$
 which on the birational model given by the Legendre family act as
 $$ y  \longleftrightarrow  -  y,  \ \ \xi  \longleftrightarrow  - \xi, \ \ \xi  \longleftrightarrow  \frac{a}{\xi}\,. $$
 
 We want now to finish explaining the monodromy of these normal forms. 
 To this purpose we consider the subgroup
 $$  \Ga_{2,4} : = \left\{ \left( \begin{array}{cc} \alpha & \beta \\ \gamma & \delta \end{array}  \right)   
\in \PP SL (2 , \ZZ)  \biggm|  \begin{array}{cc}
  \alpha \equiv 1 \mod 4, & \beta \equiv  0 \mod 4, \\\gamma \equiv  0 \mod 2, &\delta \equiv  1 \mod 2
  \end{array}  
\right\} $$
a subgroup of index $2$ of the congruence subgroup  
 $$  \Ga_2 : = \left\{ \left( \begin{array}{cc} \alpha & \beta \\ \gamma & \delta \end{array}  \right)   
 \in \PP SL (2 , \ZZ) \biggm| \begin{array}{cc}
  \alpha \equiv 1 \mod 2, & \beta \equiv  0 \mod 2, \\\gamma \equiv  0 \mod 2, &\delta \equiv  1 \mod 2
  \end{array}  
\right\} .$$

To the chain of inclusions $$ \Ga_{2,4}  <    \Ga_2  < \PP SL (2 , \ZZ) $$
corresponds a chain of fields of invariants 
$$ \CC (j) \subset \CC(\la) = \CC(\tau)^{\Ga_2} \subset  \CC(\tau)^{\Ga_{2,4}}\,, $$
where the respective degrees of the extensions are 6, 2.

Here, $\la$ is, as previously mentioned, the cross-ratio of the four points 
$\mathfrak p (0), \mathfrak p (\frac{1}{2}), \mathfrak p (\frac{\tau}{2}), \mathfrak p (\frac{1 + \tau}{2}), $
 where $\mathfrak p $ is the Weierstra\ss \  function, and $j (\la) = \frac{4}{27} \frac{(\la^2 - \la + 1)^3}{\la^2 (\la-1)^2}$ is the $j$-invariant.

If $\la (a)$ is the cross ratio of the four points $1,-1, a, - a$,    then $\la (a) = \frac{(a-1)^2}{( a+1)^2 }$, hence $ a = \frac{1}{\la -1} ( -1 - \la \pm 2 \sqrt \la)$
and  $\CC(a) = \CC(\sqrt \la)$ is a quadratic extension. 
The geometric meaning 
is related  to the algebraic formulae that we have illustrated in greater generality,  which concretely 
explain the degree 2 field extension $\sqrt \la$ as follows:
on the elliptic curve in Legendre normal form $y^2 = (x^2-1) (x^2 - a^2)$ we have the automorphism
$  g_1 (x,y) = (-x, -y)$; and   
the  quotient by $g_1$, where we set $u : = x^2 ,  \ v : = xy, \la = a^2$ is  the elliptic curve of equation
$ v^2 = u (u-1) (u - \la).$  Under this quotient map, a point of 4-torsion, $\tau/4$ on  $E_{\tau} = \CC/ \ZZ + \ZZ \tau $ is sent to a  point of 2-torsion
inside the curve $E_{\tau/2} \CC/ \ZZ + \ZZ \tau/2$.

Hence, unlike what looked to be superficially at the beginning,  the two normal forms for which  the four   2-torsion points are given an order are not the same;
and the reason behind this is that in the Legendre normal form one describes the full action of the group $(\ZZ/2)^3$, without the asymmetry of
treating the four  2-torsion points on a different footing (in the normal form with $\la$ one needs a non affine transformation in order to exchange the point
$x=0$ with the point $x=\infty$).

Setting now $b : = \sL (\frac{\tau}{4} )$, we have that 
$a = b^2$, hence $\CC(b)$ is a quadratic extension of  $\CC(\tau)^{\Ga_{2,4}}$.

In other words, the parameter $ b \in \CC \setminus \{ 0, 1, -1 , i, - i \}$  yields an unramified covering of degree $4$ of $ \la  \in \CC \setminus \{ 0, 1\}$ ,
hence the field $\CC(b)$ is the invariant field for a subgroup {{} $\Ga_{2,8}$}   
 of index $2$ in $\Ga_{2,4}$.

{{}
\noindent
By \cite[\S 182]{bianchi}, $b$ is invariant under the subgroup of $\Ga_2$ given by the transformation
such that $\alpha^2+\alpha\beta \equiv 1 \mod 8$.
Since $\alpha \equiv 1 \mod 2$, this equation is equivalent to require that $\beta \equiv 0 \mod 8$, i.e.:
$$  \Ga_{2,8} : = \left\{ \left( \begin{array}{cc} \alpha & \beta \\ \gamma & \delta \end{array}  \right)   
\in \PP SL (2 , \ZZ)  \biggm|  \begin{array}{cc}
  \alpha \equiv 1 \mod 4, & \beta \equiv  0 \mod 8, \\\gamma \equiv  0 \mod 2, &\delta \equiv  1 \mod 2
  \end{array}  
\right\} .$$
}

 In other terms,  the  family 
 $ E(b)$  is the family of elliptic curves with a   $\Ga_{2,8}$ level  structure:
 it  is  the quotient of $(\CC\times \HH)$,
 with coordinates $(z, \tau))$, by the action of the group (a semidirect product) generated by  
 $(\ZZ^2)$ which acts  by 
 $$(m, n) \circ (z, \tau) = ((z + m + n \tau, \tau)$$
 and by $ {{} \Ga_{2,8}}  \subset  \PP SL (2 , \ZZ)$.

 We have a family of  ramified covers: 
 $$ (*) \ \ \CC \setminus \{ 0, 1, -1 , i, - i \}  \ra  \CC \setminus \{ 0, 1, -1  \} \ra \CC \setminus \{ 0, 1  \} \ra \CC,$$ 
 $  \CC \setminus \{ 0, 1, -1 , i, - i \}$ with coordinate $b$ maps to $ a = b^2 \in  \CC \setminus \{ 0, 1, -1  \}$,
 in turn we have another degree two map  $$ a \in \CC \setminus \{ 0, 1, -1  \} \mapsto   \la  = \frac{(a-1)^2}{( a+1)^2 } \CC \setminus \{ 0, 1  \},$$
 finally  we have a degree 6 map $ j = \frac{4}{27} \frac{(\la^2 - \la + 1)^3}{\la^2 (\la-1)^2} $.
 
 Denote by $\HH$ the upper half plane $\HH : = \{ \tau | Im (\tau) > 0 \}$: identifying two such open sets of the complex line (as in $(*)$) as quotients 
 $\HH / \Ga$,   $\HH / \Ga'$,  with $\Ga \subset \Ga'$,  the monodromy of  $\HH / \Ga \ra \HH / \Ga'$ is the image of $\Ga'$ in the group of perutations
 of the set $\Ga' / \Ga$ of cosets of $\Ga$ in $\Ga'$. Each such map in the diagram above  is a Galois covering (but not all their compositions!).
 
 In particular, the monodromy of $ \la \ra j$ is the symmetric group $ \mathfrak S_3$ in three letters  permuting the three points $0,1, \la$: it can be seen as the group
 of linear automorphisms of the subgroup $E[2] \cong (\ZZ/2)^2$ of 2-torsion points of $E$, obtaining classical isomorphisms
  $$ \mathfrak S_3 \cong GL (2, \ZZ/2) = \PP SL (2, \ZZ) / \Ga_2.$$
  
  This is the perhaps the reason why the normal form  $$ y^2 = x (x-1)(x - \la),  \la \neq 0,1,$$
  the normal form for elliptic curves given together with an isomorphism of the group of 2-torsion points with $(\ZZ/2)^2$,
 is the most used nowadays for theoretical purposes.
 
  There is in fact no normal form which depends only upon the variable $j$,
 since, varying $j$, one would obtain a family of curves $C_j$ which are isomorphic to  the elliptic curve $E_j$ when $j\neq 0,1$,
 but for $j=0$ (and similarly for $j=1$) one would get 
 $C_0 = \PP^1$, the quotient of the elliptic curve $E_0$ with invariant $j=0$  by the cyclic group of order 3 of automorphisms acting by  
 $ z \mapsto \eta z$, with $\eta^3 =1$.
 
 Moreover, this form is useful in arithmetic because of the method of 2-descent. This is related  to the procedure, that we have just described algebraically,
 of being able to divide the elliptic curve
 by the points of 2-torsion. In fact, the   affine curve $C$
 with equations
$$  w^2 = (x^2 -1) ,  t^2 = (x^2 - a^2)$$
admits a group of automorphisms given by the translations by points of 2-torsion,
and generated by the two automorphisms
$$ h_1 (x,w,t) = (x, -w, -t) , \ h_2  (x,w,t) = (-x, w, -t).$$ 
Since the invariants are $v : = x w t , u : = x^2$, its quotient is the  curve
$$ v^2 = u (u-1) (u-a^2),$$
which is isomorphic to $C$, whereas the quotient by $h_1$ is, setting $ y : = wt$,
the elliptic curve in Legendre normal form
$$ y ^2 = (x^2 -1) (x^2 -a^2).$$
  More general algebraic formulae for torsion points of higher order were found by Bianchi (see \cite{bianchi}).

\section{Periodic functions and Abelian varieties}

Writes Reinhold Remmert on page 335 of his book \cite{remmert}, in the section dedicated to short biographies 
of the creators of function theory, in alphabetical order  Abel, Cauchy, Eisenstein, Euler, Riemann and Weierstra\ss \ :

`Karl Theodor Wilhelm Weierstra\ss \ ,  born in 1815 in Ostenfelde, Kreis Warendorf, Westphalia...studied Mathematics
in the years 1839-40 at the Academy of M\"unster, passed the State exam with Gudermann \footnote{Weierstra\ss \  attended Gudermann's lecture on elliptic functions, some of the first lectures on this topic to be given(\cite{mac})};  1842-1848,
teacher at the Progymnasium in Deutsch-Krone, West Prussia, for Mathematics, calligraphy and gymnastic; 
1848-1855, teacher at the Gymnasium in Braunsberg, East Prussia; 1854, publication of the ground breaking 
results he had  gotten already in 1849, in the article `Zur theorie der Abelschen Functionen' in vol. 47 of Crelle's Journal,
hence he receives a honorary PhD from  the University of K\"onigsberg.. ;1864, ordinarius (full professor) at the Berlin University...
1873-74, Rector of the Berlin University ..he died in 1897 in Berlin.' 

In fact, the groundbreaking results appeared in two articles on Crelle, \cite{weieab1} and \cite{weieab2}, the first devoted 
to the general representation of abelian functions as convergent power series, and the second giving a full-fledged theory of the inversion of hyperelliptic integrals.

What are the Abelian functions? They are just the meromorphic functions on $\CC^n$ whose group of periods
$$ \Ga_f : = \{ \ga  \in \CC^n | f (z + \ga) = f (z) , \ \forall z \in \CC^n \}$$
is a discrete subgroup of maximal rank $= 2n$. In other words, such that $ \Ga_f $ consists of the $\ZZ$-linear combinations of $2n$ vectors 
which form an $\RR$-basis of $\CC^n$  as a real vector space.

In general $ \Ga_f $ is a closed subgroup, which can be uniquely written as the sum of a $\CC$-vector subspace $V_1$
of $\CC^n$ with a discrete subgroup $\Ga_2$ of a supplementary subspace $V_2$ (i.e., $\CC^n \cong V_1 \oplus V_2$). 
If $V_1$ is $\neq 0$, then the function $f$ is degenerate, i.e., it depends upon fewer variables than $n$.
A central question of the theory of periodic meromorphic functions has been: given a discrete subgroup
$\Ga \subset \CC^n$, when does there exist a nondegenerate function $f$ which is $\Ga$-periodic? 

In the case of $n=1$, the answer is easy, the abelian functions are just  the elliptic functions, and, fixed a discrete subgroup of
rank $=2$,  the field of
$\Ga$-periodic  meromorphic functions is the field of meromorphic functions
on the elliptic curve $E = \CC / \Ga$, generated by $\wp $ and $\wp '$, which satisfy the unique relation
$$ (\wp ') ^2 = 4 \wp ^3 - g_2 \wp  - g_3.$$
So far, so good.

For $n\geq 2$, things are more complicated, and it took more than the ideas of Riemann and Weierstra\ss \ 
in this long and fascinating story!

An important contribution was given by the French school, for instance a basic theorem was the theorem of Appell-Humbert, 
also called the theorem on the linearization of the system of exponents, and 
proven by these authors for $n \leq 2$, and by Conforto in arbitrary dimension (\cite{conforto}, see also \cite{siegel}, page 53,
for an account of the history). 

The starting point\footnote{ writes Siegel (page 23, \cite{siegel}): this question leads to interesting and complicated problems, problems whose solutions were initiated by Weierstra\ss \  and Poincar\"e, continued by Cousin, and completed mainly by Oka.. }
is a theorem due to Poincar\'e, showing that every meromorphic function $f$ on $\CC^n$ can be written as the
quotient of two relatively prime holomorphic functions:
$$ f (z) = \frac{F_1(z)}{F_0(z)}.$$

The $\Ga$-periodicity of $f$, in view of unique factorization,  translates into the same functional equation for numerator and denominator:
$$ F_j ( z +\ga) = k_{\ga} (z) F_j (z),  $$
where the holomorphic functions  $ k_{\ga} (z)$ are nowhere vanishing, and satisfy the cocycle condition
$$   k_{\ga_1 + \ga_2} (z) =  k_{\ga_1} (z + \ga_2)  k_{\ga_2} (z).$$ 
 
 Since  $ k_{\ga} (z)$ is nowhere vanishing, one can take its logarithm, and write
 $ k_{\ga} (z) = exp (\phi_{\ga} (z))$, and the  functions $(\phi_{\ga} (z))$ are called a system of exponents.
 
 Since however the above factorization is unique only up to multiplying with a nowhere vanishing function $ exp (g(z)$,
 we see that we can replace the system of exponents  by a cohomologous one
 $$ \varphi_{\ga} (z): =  \phi_{\ga} (z) - g( z + \ga) + g(z).$$
 
 The main result of Appell-Humbert and Conforto is that one can find an appropriate function $g(z)$ such that 
 the resulting exponent $ \varphi_{\ga} (z)$ is a polynomial of degree at most one in $z$; from this it follows, by pure bilinear algebra arguments,
 that the system can be put into a unique normal form, the so called
 
 {\bf Appell-Humbert normal form for the system of exponents:}
$$  \hat{k}_{\ga} (z) = exp (\varphi_{\ga} (z))  =  \rho (\ga) \cdot exp (\pi H(z,\ga) + \frac{1}{2} \pi H(\ga,\ga) ),$$
where $H(z,w)$ is a Hermitian form satisfying the
 
{\bf First Riemann bilinear relation: the imaginary part of $H$ yields an alternating  bilinear form with integer values on $\Ga$,}
$$ \sE: = Im (H) : \Ga \times \Ga  \ra \ZZ,$$
and $\rho  : \Ga \ra S^1 : = \{ w \in \CC | |w| = 1\}$ is a semicharacter for $\sE$, i.e., we have 
$$  \rho (\ga_1 + \ga_2) =  exp ( 2 \pi i \sE(\ga_1, \ga_2))  \rho (\ga_1 )   \rho ( \ga_2) .$$

Moreover, the functional equation 
$$ (FE): \  F_j ( z +\ga) =  \hat{k}_{\ga} (z) F_j (z)  $$
has some nondegenerate solution if and only if the 

{\bf Second  Riemann bilinear relation is satisfied: $ H > 0$,  the Hermitian form  $H$ is positive definite on $\CC^n$.}

The Riemann bilinear relations allow then, after a change of variables, to assume that the group of periods
can be put in a semi-normal form, similar to the one of elliptic curves:
$$ \Ga =  T \ZZ^n \oplus  \tau \ZZ^n,$$
here $T, \tau$ are $ n \times n$ matrices, $T = diag (t_1, t_2, \dots t_n), t_i \in \NN, t_1 | t_2 | \dots | t_n$
(the $t_j$'s  are the elementary divisors of the matrix $\sE$), 
while $\tau$ is a matrix in the Siegel generalized upper half space
$$ \tau \in \HH_g : = \{ \tau \in Mat(n,n, \CC) | \tau = ^t \tau, Im (\tau) > 0 \}.  $$

This change of variables is crucial, since it allows to derive, after some algebraic manipulation,  all the solutions of the functional equation $ (FE)$
from  the Riemann theta function;  indeed the solutions are linear combinations
of  the theta functions with characteristics $a \in \QQ^n$:
\footnote{the Riemann theta function is just the one with characteristic $a=0$}
$$  \theta [a,0] (z, \tau) = \Sigma_{p \in \ZZ^n} exp (  \pi i (   ^t(p+a) \tau (p+a) + 2  ^t z (p+a))).$$ 

Nowadays, see for instance the book by Mumford (\cite{abvar}), the theorem of Appell-Humbert is viewed as the consequence
of the exponential exact sequence on the complex torus $ X : = \CC^n / \Ga$:
$$ 0 \ra H^1 (X, \ZZ) \ra H^1 (X, \hol_X) \ra H^1 (X, \hol^*_X)  \ra  H^2 (X, \ZZ) \ra H^2 (X, \hol_X) \ra ..$$
Here  $ H^2 (X, \ZZ)$ is the space of alternating forms $\sE$ as above,  the Hermitian form $ H(z,w)$ 
is the unique representative with constant coefficients of  the first Chern form   of the line bundle $L$ 
which is associated to the cocycle $k_{\ga}(z)$
as follows:
$L$ is defined as the quotient of $\CC^n \times \CC$ by the action of $\Ga$ such that 
$$ \ga (z,w) = ( z + \ga, k_{\ga}(z) w). $$ 
Finally, $H^1 (X, \hol^*_X) = : Pic(X)$ is the group of isomorphism classes of line bundles on $X$.

These new formulations through K\"ahler geometry, and Hodge theory, allowed a very vast generalization,
the beautiful theorem by Kodaira \cite{kodaira}.

{\bf  Kodaira's embedding theorem:}

{\em A compact complex manifold $X$ is projective if and only if it has a positive line bundle, i.e., a line bundle $L$ admitting a Chern form
which is everywhere strictly positive definite.}

However beautiful Kodaira's theorem, it did not completely solve the problem about the existence of nondegenerate $\Ga$-periodic
meromorphic functions, which had been dealt by Cousin (\cite{cousin}) for the case $n=2$, a case which can be reduced,
in the case where $ X : = \CC^2 / \Ga$ is noncompact,  to the theory of elliptic functions. Indeed, Cousin and Malgrange 
(\cite{malgrange}) showed
that the linearization for the system of exponents does not hold when $ n \geq 2$ and we are  in the noncompact case (where harmonic theory
can no longer be used).

Nevertheless, in \cite{cc}, we were able to show that the Riemann bilinear relations are always necessary and sufficient conditions 
on $\Ga$ for the existence of nondegenerate $\Ga$-periodic
meromorphic functions. The lack of linearization of the system of exponents, however, make the description of the $\Ga$-periodic functions
not so simple as in the compact case.

This said, we have not yet explained what is the inversion of abelian integrals. We shall do it in modern terminology.

Let $C$ be a compact Riemann surface of genus $g$, so that the space  $H^0(\Omega^1_C)$ of holomorphic one-forms 
has complex dimension equal to $g$, and $C$ is obtained  topologically from a polygon with $4g$ sides
$$\al_1,\be_1, \al_1^{-1},\be_1^{-1},\dots , \al_g,\be_g, \al_g^{-1},\be_g^{-1}$$ glueing each side $\al_j$ with the side $\al_j^{-1}$,
which is oriented in the opposite direction with respect to one chosen orientation of the boundary (and doing similarly 
for $\be_j$ and  $\be_j^{-1}$). 

{\bf The Jacobian variety $J(C)$ } is nothing else than the quotient of the dual vector space $H^0(\Omega^1_C)^{\vee}$
by the discrete subgroup $\Ga$ generated by the integrals over the closed paths in $C$,  $\al_j, \be_j$.
For an appropriate choice of the polygonal dissection and choice of a basis $\omega_1, \dots \omega_g \in H^0(\Omega^1_C)$
one obtains:
$$ \int_{\al_i} \omega_j = \de_{i,j}, \   \int_{\be_i} \omega_j  = \tau_{i,j}, \ (\exists \tau \in \HH_g).$$ 

In fact  the Riemann bilinear relations were discovered by Riemann exactly in the case of the period group of a
compact Riemann surface (and then extended to the general case).

The Abel-Jacobi  map (later generalized by Francesco Severi and Giacomo Albanese to the case of higher dimensional varieties) is the map
$$ a : C \ra J(C),  a (P) : = \int_{P_0}^P \in  H^0(\Omega^1_C)^{\vee} / \Ga .$$ 

This map can be extended to the Cartesian product $C^m$ by $$ a (P_1, \dots, P_m) : = \Sigma_{j=1}^m a (P_j)$$
and clearly  its value is independent of the order of the points $P_j$, hence, once we define the m-fold symmetric product  
$ C^{(m)}$ as the quotient $C^m / \mathfrak S_m$, we have a similar  Abel-Jacobi  map $ a_m  : C^{(m)} \ra J(C)$. In these terms
one can express the

{\bf Jacobi inversion Theorem: $ a_g  : C^{(g)} \ra J(C)$ is surjective and 1-1 on a dense open set.}

As we already mentioned, Jacobi did not prove this theorem, just posed the problem in 1832,  and Weierstra\ss \  proved it in the special case where
$C$ is a hyperelliptic curve, i.e.,  $C$ is a double cover of the projective line, i.e., given by an equation
$$ y^2 = \Pi_{j=1}^{2g+2} ( x - \la_j) ,$$ 
quite reminiscent of the standard equation for elliptic curves.

If however one looks at the older literature, statements are harder to follow because there the authors  speak of 
abelian functions, abelian integrals, whereas nowadays the statement is made in terms of varieties, and in the language of algebraic geometry.

Hence our question: where and when did this geometrization process begin?

Certainly the French school gave a big impetus: they used the term `Hyperelliptic variety' for a variety which is the quotient 
of $\CC^n$ by a group $\Ga$ acting freely, but not necessarily a group of translations. Whereas the name `Picard variety' 
was reserved for  what we call now Abelian varieties (except that the Picard variety of an Abelian variety $A$  is the dual 
of $A$!). The reader can find more details in the article by Kleiman \cite{kleiman}.

But, for the inversion of integrals, a problem which  Appell and Goursat in their book (\cite{a-g}, page 463) 
acknowledge as solved by Riemann and Weierstrass simultaneously, comes into play the Jacobian variety of a curve for a geometric formulation:
the Jacobial variety $J(C)$ of a curve of genus $g$ is birational to the g-fold symmetric product of the curve $C$. 

The term `Variet\'a di Jacobi' appears clearly in the title
of the famous 1913 article by Ruggiero Torelli \cite{torelli}, where he proved the famous

{\bf Torelli's theorem:  two curves are  isomorphic if and only if their  Jacobian varieties are isomorphic as polarized Abelian varieties.}

Indeed, the term had been used a few years earlier by Enriques and Severi (\cite{e-s}), who spoke of the Jacobian surface of a curve 
of genus 2 while describing the classification of hyperelliptic surfaces. Mathematically, the inversion theorem had been proven, amidst other results,  by Castelnuovo in \cite{castelnuovo} as a consequence
of the Riemann-Roch theorem for curves. However, at that time Castelnuovo did not yet speak of the `Jacobian variety';
rather, when his selected papers were reprinted in 1936, he added a note at the end of the paper:
`the results of the paper can be stated in the simplest way introducing the concept of the g-dimensional Jacobian variety of the curve $C$..'.

Research on Abelian varieties and on monodromy is still an active ongoing field (see \cite{kaw}, \cite{del}, \cite{DeligneMostow},..).

\subsection{Cohomology as difference equations}
 It is difficult to give a flavour of the many intriguing questions which are involved in the development of the theory of Abelian functions and varieties,
 so I shall try with an elementary example, which is equivalent to the vanishing of the cohomology group $H^1 (\CC^*, \hol_{\CC^*})=0$,
 and which cannot be superseded by harmonic theory when $\Ga$ is not of maximal rank.

 \begin{lemma}
 Let $f(z)$  be an entire  holomorphic function on $\CC$: then there exists another entire holomorphic function $g(z)$ on $\CC$ such that
 $$    f(z)  = g(z+1) - g(z). $$ 
 \end{lemma}
 
 At first glance, it seems like that the correct approach would be  to use that entire functions are just sums of a power series,
  converging everywhere on $\CC$.
  
  That is, we may write (keeping track that a constant  for $g$ yields the result $0$):
  $$ f(z) = \Sigma_{i=0}^{\infty} b_i z^i , \  \  g(z) = \Sigma_{j=1}^{\infty} a_j z^j .$$
  
  The inhomogeneous equations to solve are given by
  $$ b_i =  \Sigma_{j=1}^{\infty}  {{j}\choose{i}} a_j ,$$
   and are unfortunately given by an infinite sum.
   
   However, if we restrict to the case where $f$ is a polynomial of degree  $r$, then there exists a unique polynomial $g$ of degree $r+1$,
   and with vanishing constant term $g(0)$, which yields the solution: since we have an upper triangular invertible matrix.
   
   The  hope that the unique solutions $g_{h+1}$ for the partial sums $f_h$ converge to a  unique power series solution $g$
   is discarded once we observe that the associated homogeneous equation $g(z+1) - g(z)= 0$ has an infinitely dimensional space of solutions,
   given by all the holomorphic Fourier series $\Sigma_{h \in \ZZ} c_h w^h$, where we set $ w = exp (2 \pi i z)$.
   
   However, (see \cite{siegel}, Theorem 2, pages 49-53) a very clever use of the Schwarz lemma \footnote{which, citing again Caratheodory's preface to \cite{caratheodory},`
   allows for new type of arguments such as were unknown prior to its discovery'} and of the Weierstra\ss \  infinite products
   permits to reduce to the case where $f$ is a polynomial, hence the solution can be found.
   
   There is a second proof, which is most interesting, and based on  the use of complex analysis to solve the difference equation
   (I am illustrating it since these methods played a pivotal role in  complex analysis during  the 20th century) , as follows.
   
   Step I. One can find first a $\sC^{\infty}$-solution $g$, observing that indeed it suffices to find the solution
   for $ z \in \RR$ (since,  if we write $ z = x + i y$, then $z+1 = (x+1) + i y$).
   This is done observing that we may first determine $g$ on the interval $ [-\frac{1}{2}, \frac{1}{2}]$ choosing $g$ identically zero in a neighbourhood 
   of $-\frac{1}{2}$, and equal to $ f(x-1)$  identically  in a neighbourhood 
   of $\frac{1}{2}$. Then we extend the definition of $g$ on the other intervals $ [-\frac{1}{2}+ m, \frac{1}{2}+ m]$ by using the derived formulae:
   $$  g(z + m) = g(z) + f(z) + f (z+1 ) + \dots + f (z + m-1).$$ 
   
Step II. Since $f$ is holomorphic, $u  : = \bar{\partial} g(z) = \bar{\partial} g(z+1) = u (z+1)$, hence $u$ is the pull-back of a 
$\sC^{\infty}$ function $U$ on $\CC^*$.
It suffices to show that the equation $ U = \bar{\partial} \Psi$ is solvable in $\CC^*$, since then $ \psi (z) : = \Psi (exp(2 \pi i z))$
satisfies 
$$   \bar{\partial} g(z) =  \bar{\partial} \psi (z), \ \psi (z+1) = \psi (z), $$ 
hence $g_2 : = g - \psi$ is the desired holomorphic solution ($ \bar{\partial} g_2(z) \equiv 0$) .

Step III. The  equation $ U = \bar{\partial} \Psi$ is solvable in $\CC^*$ when $U= U_j$ has compact support, say contained 
in $K_j = \{ z | 1/j \leq |z| \leq j\}$, where $j$ is an integer. It suffices to use the Bochner-Severi-Martinelli integral formula
( see H\"ormander, \cite{hormander} theorem 1.2.2, page 3)
$$  \psi_j (z) = \frac{1}{ 2 \pi i} \int \frac{u_j(\zeta)} { (\zeta - z) } d \zeta d \bar{\zeta}.$$ 

Step IV.  Let $\phi_j(z)$ be a function which is $\sC^{\infty}$, identically $\equiv 1$ in a neighbourhood of $K_j$, and with support $\subset K_{j+1}$,
and set $v_j := \phi_j(z) - \phi_{j-1}(z)$.

Then by step III there is a function $\psi_j$ such that $\bar{\partial}  \psi_j = U_j : = v_j U$. Since $v_j \equiv 0$ on $K_{j-1}$,
$\psi_j$ is then holomorphic in a a neighbourhood of $K_{j-1}$, hence there exists a global Laurent series $\tilde{\psi_j }$
which approximates $\psi_j$ , say that
$$  | \tilde{\psi_j } -\psi_j | < 2^{-j} \ on \  K_{j-1}.$$ 

We finally define $\psi : = \Sigma_j  ( \psi_j  - \tilde{\psi_j } )$, where the series is convergent in norm over compact sets.
In particular, the sum $\Sigma_{j \geq h+1} ( \psi_j  - \tilde{\psi_j } )$ is a holomorphic function in a neighbourhood of $K_h$.

So, $\psi $ is $\sC^{\infty}$ and 
$$ \bar{\partial}  \psi =   \Sigma_j   \bar{\partial}    ( \psi_j  - \tilde{\psi_j } ) = \Sigma_j {U_j } =  U ( \Sigma_j v_j )  = U.$$ 

\bigskip
\noindent
{\bf Acknowledgements.}
I would like to thank Frans Oort for raising my interest concerning the history of the monodromy theorem,
and introducing me to some psychological  aspects of  mathematics. And other friends who explained to me
that the difference between mathematics and history of mathematics  is the same as between everyday's life
and a court case: in court facts, and not ideas, count, and for everything one needs evidence (really so?).



\begin{thebibliography}{Grif-SchmX}

\bibitem[A-G895]{a-g}
Appell, Paul; Goursat, Edouard,
{\em Th\'eorie des Fonctions alg\'ebriques et de leurs int\'egrales. \'Etude des fonctions analytiques sur une surface de Riemann.}
Gauthier-Villars et Fils, Paris (1895), X + 530 pages.

\bibitem[BC11a]{BC11burniat1}
Bauer, Ingrid; Catanese, Fabrizio,
{\em  Burniat surfaces {I}: fundamental groups and moduli of primary
  {B}urniat surfaces},
in {\em Classification of algebraic varieties}, EMS Ser. Congr. Rep.,
  pages 49--76. Eur. Math. Soc., Z\"urich (2011).
  
  \bibitem[BCF15]{bcf}
  Bauer, Ingrid; Catanese, Fabrizio;   Frapporti, Davide,  
  {\em Generalized Burniat type surfaces and Bagnera-de Franchis varieties}
  J. Math. Sci. Univ. Tokyo 22,  Kodaira Centennial issue, 1--57 (2015). 
  
  
\bibitem[Bia01]{bianchi}
Bianchi, Luigi,
{\em Lezioni sulla teoria delle funzioni di variabile complessa e delle funzioni ellittiche,} Enrico Spoerri, Pisa, 
 630 pages, (1901); second edition, 685 pages (1916).
 
 \bibitem[CapCat91]{cc}
 Capocasa, Francesco; Catanese, Fabrizio,
 {\em Periodic meromorphic functions,}
 Acta Math., 166 (1991), 27-68.
 
\bibitem[Car50]{caratheodory}
C. Caratheodory,
{\em Funktionentheorie}, 2 Volumes, Birkh\"auser,  (1950); english translation: Theory of functions of a complex variable, 2 vol., Chelsea Publ.,
xii + 304, resp. 220 pages (1954)

\bibitem[Cart]{cartan}
Cartan Henri,
{\em Th\'eorie el\'ementaire des functions analytiques d' une ou plusieurs variables complexes}, 4th edition, Hermann Paris, 231 pages (1961).

\bibitem[Cast93]{castelnuovo}
Castelnuovo, Guido,
{\em  Le corrispondenze univoche tra gruppi di $p$ punti sopra una curva di genere $p$, }
Rendiconti del R. Istituto Lombardo, s. $2^a$, vol. 25, (1893); reprinted (pages 79-94) in:
Guido Castelnuovo, {\em Memorie scelte,} with notes by Castelnuovo, Zanichelli, Bologna x + 588 (1937).

\bibitem[Con42]{conforto}
Conforto, Fabio,
{\em Funzioni abeliane e matrici di Rieman, Parte I} Corsi Ist. Naz. Alta Matematica, Universit\'a di Roma (1942);
later version {\em Abelsche Funktionen und algebraische Geometrie,}
Springer-Verlag, Berlin-G\"ottingen-Heidelberg (1956).

\bibitem[Cou910]{cousin}
Cousin, P.
{\em Sur les functions triplement periodiques de deux variables,}
Acta Math. 33, 105-232 (1910).

\bibitem[Del70]{del}
P. Deligne,
{\em \'Equations diff\'erentielles \'a points singuliers r\'eguliers,}
 Lecture
Notes in Mathematics, Vol. 163. Springer-Verlag, Berlin-New 
York (1970),  iii+133 pp.


\bibitem[D-M86]{DeligneMostow}
P.~Deligne and G.D. Mostow.
\newblock Monodromy of hypergeometric functions and non-lattice integral
  monodromy.
\newblock {\em Publ. Math. IHES}, 63:\ 5--89, (1986).

\bibitem[Ency09]{wirtinger}
H. Burkhardt, W. Wirtinger, R. Fricke, E. Hilb
{\em Encyklop\"adie der Mathematischen Wissenschaften mit Einschluss ihrer Anwendungen,} vol. 2, Analysis,
Vieweg+Teubner Verlag, 674 pages, (1909).

\bibitem[Eis]{eisenstein}
G. Eisenstein,
{\em Mathematische Abhandlungen,} Georg Olms Verlagsbuchhandlung, Hildesheim (1967).

\bibitem[E-S907]{e-s}
Enriques, Federigo; Severi, Francesco,
{\em Intorno alle superficie iperellittiche, }
Rend.  Acc.  Lincei XVI (1907), 443-453.

\bibitem[Frei06]{freitag}
Freitag,Eberhard, Busam Rolf,
{\em Funktioentheorie 1}, Springer Verlag, (2006).

\bibitem[F-K912]{f-k}
Fricke, Robert; Klein, Felix,
{\em Vorlesungen \"uber die Theorie der Automorphen Functionen,} 
vol. 2 {\em Die Functionentheoretischen Ausf\"uhrungen und die Anwendungen},
Verlag von B.G. Teubner, Leipzig und Berlin, viii + 668 pages (1912).

\bibitem[God58]{godement}
Godement, Roger {\em Topologie alg\'ebrique et th\'eorie des faisceaux. }  Actualit\'es Scientifiques et Industrielles 
No. 1252. Publ. Math. Univ. Strasbourg. No. 13 Hermann, Paris (1958) viii+283 pp. ; reprinted (1973).

\bibitem[H\"orm73]{hormander}
H\"ormander, Lars,
{\em An introduction to Complex Analysis in Several Variables,}
North-Holland Mathematical Library, vol. 7 , x + 213 pages, North-Holland, Amsterdam, London (1973).

\bibitem[Hur83]{hur-rat}
Hurwitz, A.,
{\em Proof of the theorem that a single valued function of arbitrarily many variables that is everywhere representable as a quotient of two power series is a rational function of its arguments. (Beweis des Satzes, dass eine einwertige Function beliebig vieler Variabeln, welche \"uberall als Quotient zweier Potenzreihen dargestellt werden kann, eine rationale Function ihrer Argumente ist.)} 
Kronecker J. XCV, 201--207 (1883).

\bibitem[Hur19]{hurwitz}
Hurwitz, Adolf,
{\em Vorlesungen \"uber allgemeine Funktionentheorie und elliptische Funktionen (Lectures on general function theory and elliptic functions),}
5th edition,( With an introduction by  Reinhold Remmert)  Berlin: Springer. xxiv, 249 S.  (2000).

\bibitem[Hur-C22]{h-c}
Hurwitz, Adolf,
{\em Vorlesungen \"uber allgemeine Funktionentheorie und elliptische Funktionen, } edited and with the addition of a section on {\em geometrische Funktionentheorie}
by  Richard Courant.. 
Berlin: J. Springer, XI u. 399 pages (1922).


\bibitem[Hur91]{hurmoves}
Hurwitz, A.,
{\em On Riemann surfaces with given branch points. (Ueber Riemannsche Fl\"achen mit gegebenen Verzweigungspunkten.)}
Math. Ann. XXXIX, 1--61 (1891).

\bibitem[Ino94]{inoue}
M.~Inoue.
\newblock Some new surfaces of general type.
\newblock {\em Tokyo J. Math.}, 17(2):295--319, 1994.

\bibitem[Kaw81]{kaw}
Kawamata, Yujiro {\em Characterization of abelian varieties.} Compositio Math. 43 (1981), no. 2, 253--276.


\bibitem[Klei04]{kleiman}
Kleiman, Steven L.,
{\em What is Abel's theorem anyway?}, in
`The Legacy of Niels Henrik Abel',
The Abel bicentennial, 2002, Laudal and Piene editors, Springer Verlag (2004), 395--440.

\bibitem[Kod54]{kodaira}
Kodaira, Kunihiko,
{\em On K\"ahler varieties of restricted type (an intrinsic characterization of algebraic varieties),}
Annals of Math. , 60 (1954), 28--48. 

\bibitem[Kraz03]{krazer}
Krazer, Adolf, 
{\em Lehrbuch der Thetafunktionen,}
Leipzig: B. G. Teubner, x + 509 pages (1903),
Chelsea reprint (1970).


\bibitem[Lag797]{lagrange}
 Lagrange Joseph Louis,
  {\em Th\'eorie des fonctions analytiques },  Impr. de la R\'epublique, Paris, Prairial an V (1797).
Second edition: Courcier, Paris, (1813),also in Lagrange's  Oeuvres, vol. IX.
 Oeuvres de Lagrange, ed. by M. J.-A. Serret [et G. Darboux], Gauthier-Villars, Paris, (1867-
1892),14 volumes.

\bibitem[Malg75]{malgrange}
{\em La cohomologie d' une variet\'e analytique \`a bord pseudoconvexe
n' est pas necessairement separ\'ee,}
C.R. Acad. Sci. Paris, 280 (1975), 93--95.
 
 \bibitem[Mum70]{abvar}
D.~Mumford.
\newblock {\em Abelian varieties}.
\newblock Tata Institute of Fundamental Research Studies in Mathematics, No. 5.
  Published for the Tata Institute of Fundamental Research, Bombay; Oxford
  University Press, London, 1970.
  
  \bibitem[O'C-R98]{mac}
  O' Connor, J.J., Robertson, E.F.,
  {\em Weierstra\ss \ ' biography,} Mac Tutor, School of Mathematics and Statistics,
  University of St. Andrews, Scotland (http.www-groups.dcs.st-and.ac.uk/~history/Biographies).
 
 \bibitem[Osg07]{osgood}
Osgood, W. F.,
{\em Lehrbuch der Funktionentheorie.} In 2 volumes, vol.1 
Leipzig: B. G. Teubner. XII u. 642 pages. gr. 8? (1907).


 \bibitem[Rem84]{remmert}
 Remmert, Reinhold,
 {\em Funktionentheorie I,}
 Grundwissen Mathematik 5, , Springer Verlag, Berlin Heidelberg (1984), 2-te Auflage xvi + 360 pages (1989)
 (contains many interesting historical comments).

\bibitem[
[Rud70]{rudin}
Rudin, Walter,
{\em Real and complex analysis},
McGraw-Hill series in higher mathematics, xi +  412 pages (1970).

\bibitem[Schw73]{Schwarz}
H.A. Schwarz.
{\em \"Uber diejenigen F\"alle in welchen die Gaussische
  hypergeometrische Reihe eine algebraische Funktion ihres vierten Elements
  darstellt}.
 Journal Reine u. Angew. Math. , 75 (1873), 292--335.
 
 \bibitem[FAC55]{fac}
 Serre, Jean Pierre,
 {\em Faisceaux alg\'ebriques coh\'erents},
 Annals of Math. vol. 61, 2, pp. 197--278 (1955).
 
 
 \bibitem[Sieg]{siegel}
 
Siegel, Carl Ludwig,
{\em Topics in Complex Function Theory,}
vol. III, {\em `Abelian Functions and Modular Functions of Several Variables',} translation by Gottschling and Tretkoff
of Siegel's lecture notes in G\"ottingen {\em (Vorlesungen \"uber ausgew\"ahlte Kapitel der Funktionentheorie)} (1964-65),
Interscience tracts in Pure and Applied Math., n. 25, Wiley Interscience, ix + 244 (1973).
Translation by Gottschling and Tretkoff
of Siegel's lecture notes in G\"ottingen {\em (Vorlesungen \"uber ausgew\"ahlte Kapitel der Funktionentheorie)} (1964-65),
where  many theorems of Weierstra\ss \  are explained;  contains a vast bibliography prepared by E. Gottschling.


 \bibitem[Tor13]{torelli}
 Torelli, Ruggiero,
 {\em Sulle variet\'a di Jacobi, I - II,}
 Rend. R. Acc. Lincei (V) 22-2 (1913),  98--103, 437-441.
 Reprinted in: {\em Collected papers of Ruggiero Torelli,}
 Ciliberto, Ribenboim and Sernesi editors, Queen's papers in Pure and
 Applied Mathematics, Vol. 101, Kingston Ontario Canada, xii + 224 pages (1995).
 
 \bibitem[Tric36]{tricomian}
 Tricomi, Francesco,
 {\em Funzioni analitiche,}
 Zanichelli, Bologna, vi + 110  (1936).
 
 
  \bibitem[Tric51]{tricomiell}
  Tricomi, Francesco,
 {\em Funzioni ellittiche,}
 Zanichelli, Bologna, ix + 339  (1951).


 \bibitem[Verd67]{verdier}
 Verdier, Jean-Louis,
 {\em Des  Cat\'egories  D\'eriv\'ees des Cat\'egories Ab\'eliennes},  Ast\'erisque ,  Soci\'et\'e Math\'ematique de France, Marseilles, 239 (1996).

\bibitem[Weil38]{weil}
A. Weil, 
{ \em G\'en\'eralisation des functions ab\'eliennes, }J. Math. Pures Appl. (9) 17 (1938), 47--87.

\bibitem[Weie54]{weieab1}
Weierstra\ss , Karl,
{\em Zur Theorie der Abelschen Functionen,}
J. Reine Angew. Math. 47, 289--306 (1854).

\bibitem[Weie56]{weieab2}
Weierstra\ss , Karl,
{\em Theorie der Abelschen Functionen, }
J. Reine Angew. Math. 52, 285--380 (1856).



\bibitem[Weie80]{weieperiodic}
Weierstra\ss , Karl,
{\em  Investigations on 2r times periodic functions of r variables, (Untersuchungen \"uber die 2r-fach periodischen Functionen mit r Ver\"anderlichen) }
Borchardt J. LXXXIX, 1--8 (1880).

\bibitem[Weie78]{hurwitz-w}
Weierstra\ss, Karl,
{\em Einleitung in die Theorie der analytischen Funktionen. } (Introduction to the theory of analytic functions. Lecture Berlin 1878, written down by A. Hurwitz, worked out by P. Ullrich). 
Deutsche Mathematiker-Vereinigung e. V., Freiburg im Breisgau. Dokumente zur Geschichte der Mathematik, 4. Braunschweig (FRG) etc.: Friedr. Vieweg - Sohn. xxvii + 184 pages  (1988).



\bibitem[WeieWerke1-6]{werke}
Weierstra\ss \ , Karl,
{\em Mathematische Werke,  herausgegeben unter Mitwirkung einer von der K\"oniglich Akademie der Wissenschaften eingesetzten Kommision,}
Berlin: Mayer  M\"uller.

Vol. 1, {\em Abhandlungen I}, VIII +  356 pages(1894).

Vol. 2, {\em Abhandlungen II}, VI + 363 pages (1895).

Vol. 3, {\em Abhandlungen III} VIII + 362 pages (1903).

Vol. 4, {\em Vorlesungen \"uber die Theorie der Abelschen Transzendenten, bearbeitet von G. Hettner und J. Knoblauch}, XIV + 631 pages (1902).

 Vol. 5, {\em Vorlesungen \"uber die Theorie der elliptischen Funktionen, bearbeitet von J. Knoblauch,}VIII + 327 pages (1915).

  Vol. 6, {\em Vorlesungen \"uber Anwendungen der elliptischen Funktionen, bearbeitet von R. Rothe, } IX + 355 pages (1915).
  
  \bibitem[WeieWerke7]{werkeVR}
Weierstra\ss \ , Karl,
{\em Mathematische Werke,  herausgegeben unter Mitwirkung einer von der Preu\ss ischen Akademie der Wissenschaften eingesetzten Commission,}
Vol. 7, {\em  Vorlesungen \"uber Variationsrechnung,} bearbeitet von R. Rothe. 
 Leipzig, Akademische Verlagsgesellschaft,  VI+324 pages (1927).



\bibitem[Weyl13]{weyl}
Weyl, H.,
{\em Die Idee der Riemannschen Fl\"ache. }
Leipzig: B. G. Teubner. X + 169 pages (Math. Vorles. a. d. Univ. G\"ottingen Nr. V.) (1913).
Reprinted in:
Weyl, Hermann; Remmert, Reinhold ,
{\em The idea of the Riemann surface.}
Teubner-Archiv zur Mathematik. Supplement. 5. Stuttgart: B. G. Teubner. xiv, 240 S.  (1997).
\end{thebibliography}
\end{document}